# Data-Driven Contextual-Aware Uncertainty Set for Robust Dispatch of Power Systems

Zhaojun Ruan, Yulin Liu, Le Fu, and Libao Shi, *Senior Member, IEEE*

***Abstract*--Both the level of conservativeness and the computational burden in robust optimization are critically influenced by uncertainty set design. However, contextual side information is rarely exploited in robust dispatch of power systems characterized by irregular data distributions, which hinders the explicit characterization of the relationship between covariates and uncertain parameters. To address this issue, a data-driven method for constructing contextual-aware uncertainty set is proposed in this letter. Based on a conditional Gaussian mixture model, a set of covariates is leveraged as side information to design uncertainty sets tailored to historical data exhibiting irregular distributions. The resulting set is formulated as a union-of-subsets formulation, and a mixed integer linear reformulation is adopted to describe the worst-case realization across all subsets. Finally, the effectiveness of the proposed method is demonstrated through numerical experiments applied to robust unit commitment.**

***Index Terms*—Uncertainty set, contextual robust optimization, covariates, side information, Gaussian mixture model.**

## I. Introduction

ROBUST optimization (RO) provides an effective framework for robust dispatch of power systems [1]. However, excessively conservative decisions can be caused by overly stringent protection against low-probability scenarios [2]. Consequently, the key to RO lies in the construction of an uncertainty set that balances conservativeness with computational tractability [3].

Conventional uncertainty sets, such as box, budget, and ellipsoidal sets, are typically sized through predefined hyperparameters [4]. Recently, complex dependence structures among uncertain parameters have been captured by data-driven approaches [5], [6]. Moreover, the nonconvex boundaries induced by irregular data distributions can be characterized by these methods [2], [7]. In contrast, conventional sets are often inadequate for representing such features, potentially causing critical scenarios to be overlooked. Additionally, side information is often ignored by these methods [8], resulting in overly large uncertainty sets that are difficult to adapt across varying operating conditions, thus inducing additional conservativeness [9]. By contrast, the relationship between covariates and uncertain parameters is explicitly modeled by contextual optimization [8], [9], enabling the construction of smaller and more accurate uncertainty sets conditional on observed external information. However, the application of this approach in robust dispatch of power systems remains limited.

To overcome this limitation, a data-driven contextual-aware uncertainty set (CAUS) is developed in this letter, incorporating covariate information. Specifically, based on a conditional Gaussian mixture model (GMM), the conditional distribution is represented by the uncertainty set formulated as the union-of-subsets (UoS), a formulation capable of capturing irregular data distributions. In comparison to conventional enumeration, a mixed integer linear reformulation of this uncertainty set is derived, reducing the computational complexity associated with UoS formulations from exponential to linear growth in terms of formulation size.

Z. Ruan, Y. Liu, L. Fu, and L. Shi are with the Tsinghua Shenzhen International Graduate School, Tsinghua University, Shenzhen 518055, China (e-mail: shilb@sz.tsinghua.edu.cn).

## II. Contextual Uncertainty Set

### A. Conditional GMM With Side Information

In this letter, let $\tilde{x} \in \mathbb{R}^n$ denote the side information (e.g., forecasts and other contextual side information) and $\tilde{\xi} \in \mathbb{R}^m$ denote the uncertainty vector. Their joint distribution of $(\tilde{x}, \tilde{\xi})$ is modeled utilizing a $K$-component GMM:

$$\mathbb{G}\left(\tilde{x}, \tilde{\xi}\right) := \sum_{k=1}^{K} p^k \mathcal{N}\left(\boldsymbol{\mu}^k, \boldsymbol{\Sigma}^k\right) \tag{1}$$

where $p^k$ denotes the weight. $\boldsymbol{\mu}^k = \left(\boldsymbol{\mu}_x^k, \boldsymbol{\mu}_\xi^k\right) \in \mathbb{R}^{n+m}$ denotes the mean vector of the $k$-th component of $\mathbb{G}(\tilde{x}, \tilde{\xi})$ . The covariance matrix of the $k$-th component is denoted by:

$$\boldsymbol{\Sigma}^k = \begin{bmatrix} \boldsymbol{\Sigma}_{xx}^k & \boldsymbol{\Sigma}_{x\xi}^k \\ \boldsymbol{\Sigma}_{\xi x}^k & \boldsymbol{\Sigma}_{\xi\xi}^k \end{bmatrix} \tag{2}$$

where the subscripts $x$ and $\xi$ indicate the components of the mean vector and covariance matrix associated with $\tilde{x}$ and $\tilde{\xi}$ , respectively.

Given that the random vector $\left(\tilde{x}, \tilde{\xi}\right)$ follows a GMM $\mathbb{G}(\tilde{x}, \tilde{\xi})$ , the conditional distribution of $\mathbb{G}(\tilde{\xi} \mid \tilde{x} = x)$ given $\tilde{x} = x$ also follows a GMM, owing to the conditional invariance property of the GMM [10]. This conditional distribution is expressed as:

$$\mathbb{G}\left(\tilde{\xi} \mid \tilde{x} = x\right) := \sum_{k=1}^{K} p_{\xi|x}^k \mathcal{N}\left(\boldsymbol{\mu}_{\xi|x}^k, \boldsymbol{\Sigma}_{\xi|x}^k\right) \tag{3}$$

$$p_{\xi|x}^k = p^k \frac{\mathcal{N}\left(x \mid \boldsymbol{\mu}_x^k, \boldsymbol{\Sigma}_{xx}^k\right)}{\sum_{j=1}^{K} p^j \mathcal{N}\left(x \mid \boldsymbol{\mu}_x^k, \boldsymbol{\Sigma}_{xx}^k\right)} \tag{4}$$

$$\boldsymbol{\mu}_{\xi|x}^k = \boldsymbol{\mu}_\xi^k + \boldsymbol{\Sigma}_{\xi x}^k \left(\boldsymbol{\Sigma}_{xx}^k\right)^{-1} \left(x - \boldsymbol{\mu}_x^k\right) \tag{5}$$

$$\boldsymbol{\Sigma}_{\xi|x}^k = \boldsymbol{\Sigma}_{\xi\xi}^k - \boldsymbol{\Sigma}_{\xi x}^k \left(\boldsymbol{\Sigma}_{xx}^k\right)^{-1} \boldsymbol{\Sigma}_{x\xi}^k \tag{6}$$

With the parameters of the GMM $\mathbb{G}(\tilde{\xi} \mid \tilde{x} = x)$ , a UoS uncertainty set can be constructed from the means $\boldsymbol{\mu}_{\xi|x}^k$ and covariance matrices $\boldsymbol{\Sigma}_{\xi|x}^k$ of its $K$ Gaussian components.

### B. UoS Construction of the CAUS

The uncertainty set associated with two wind farm power outputs is utilized as an illustrative example. At period $t$, the two wind farm power outputs are represented by a two-dimensional vector $\xi_t$. The box uncertainty set at period $t$ is illustrated by the blue rectangular region in Fig. 1(a) and is formulated as:

$$\mathcal{W}_{\text{Box}} = \mathcal{W}_{\text{Box}}^1 \otimes \mathcal{W}_{\text{Box}}^2 \ldots \otimes \mathcal{W}_{\text{Box}}^T, \quad \mathcal{W}_{\text{Box}}^t = \left\{ \tilde{\xi}_t \mid \underline{\xi}_t \le \tilde{\xi}_t \le \overline{\xi}_t \right\} \quad (7)$$

where $T$ denotes the total dispatch period. $\otimes$ is the Cartesian product. $\underline{\xi}_t$ and $\overline{\xi}_t$ are the lower and upper boundaries of two wind farm power outputs, respectively.

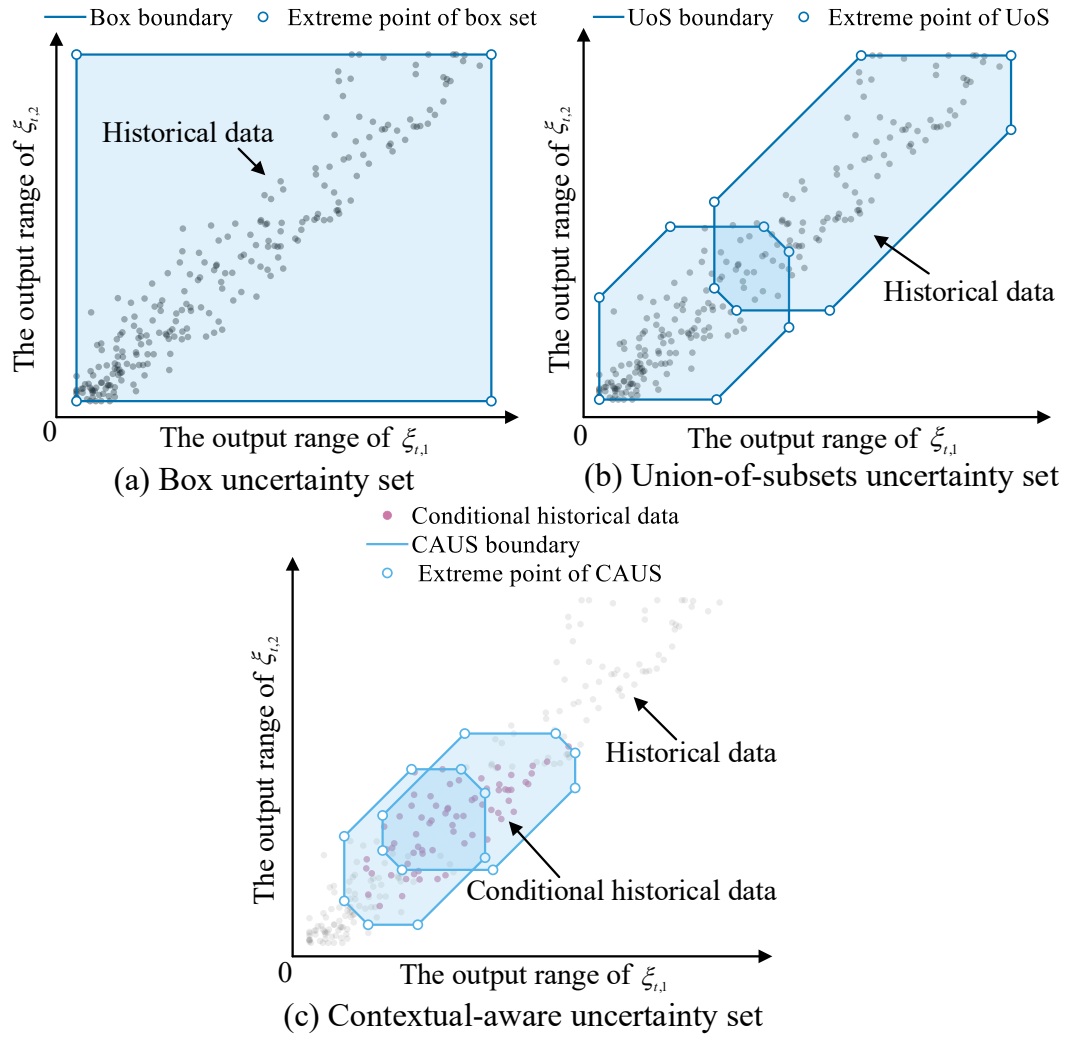

Fig. 1. Comparisons of different uncertainty sets.

To better represent irregular data distributions, the uncertainty set can instead be modeled as the UoS. As shown in Fig. 1(b), the resulting set at period $t$ is bounded by the region enclosed by two blue polygons, and the entire uncertainty set can be formulated as:

$$\mathcal{W}_{\text{UoS}} = \mathcal{W}_{\text{UoS}}^1 \otimes \mathcal{W}_{\text{UoS}}^2 \ldots \otimes \mathcal{W}_{\text{UoS}}^T, \quad \mathcal{W}_{\text{UoS}}^t = \bigcup_{k=1}^{K} \mathcal{W}_{\text{UoS}}^{k,t} \quad (8)$$

$$\mathcal{W}_{\text{UoS}}^{k,t} = \left\{ \tilde{\xi}_t \mid \tilde{\xi}_t = \boldsymbol{\mu}_k + \Lambda \Sigma^{1/2} \boldsymbol{\eta}, \|\boldsymbol{\eta}\|_\infty \le 1, \|\boldsymbol{\eta}\|_1 \le \boldsymbol{\Phi}_k \right\} \quad (9)$$

where $\boldsymbol{\eta}$ is the latent uncertainty, and the parameters $\boldsymbol{\Phi}_k$ and $\Lambda$ are utilized to control the size of the uncertainty set [7].

Although irregular historical data distributions can be described with reasonable accuracy by this method, an explicit principle for selecting the uncertainty set is not provided. Moreover, the impact of covariates or side information is not incorporated, and the resulting set selection remains relatively conservative. Based on the conditional GMM $\mathbb{G}(\tilde{\xi} \mid \tilde{x} = x)$, and given its parameters $(\boldsymbol{\mu}_{\xi|x}^k, \boldsymbol{\Sigma}_{\xi|x}^k)$, a set of $N_s$ independent and identically distributed (i.i.d.) conditional samples, denoted as $\{\xi_t^{(i)}\}_{i=1}^{N_s} \sim \mathbb{G}(\tilde{\xi}_t \mid \tilde{x} = x)$ can be drawn. For each sample, the union Mahalanobis score is defined as:

$$C_t^{(i)}(x) := \min_{k=1,2,\ldots,K} \left(\xi_t^{(i)} - \boldsymbol{\mu}_{\xi|x}^k\right)^{\mathrm{T}} \left(\Sigma_{\xi|x}^k\right)^{-1} \left(\xi_t^{(i)} - \boldsymbol{\mu}_{\xi|x}^k\right), \quad i = 1, \ldots, N_s \quad (10)$$

At period $t$, let $C_t^{(1)}(x) \le C_t^{(2)}(x) \le \ldots \le C_t^{(N_s)}(x)$ be the order statistics of $\{C_t^{(i)}(x)\}$. Given a confidence level $1-\varepsilon$, the calibrated radius parameter is set as:

$$\Gamma_{\xi|x}^t := C_t^{(\kappa)}(x), \quad \kappa := \left\lceil (1-\varepsilon)(N_s + 1) \right\rceil \quad (11)$$

The $k$-th component subset is defined by:

$$\hat{\mathcal{W}}_{\text{CAUS}}^{k,t}(x) = \left\{ \tilde{\xi}_t \in \mathbb{R}^m \mid \left(\tilde{\xi}_t - \boldsymbol{\mu}_{\xi|x}^k\right)^{\mathrm{T}} \left(\Sigma_{\xi|x}^k\right)^{-1} \left(\tilde{\xi}_t - \boldsymbol{\mu}_{\xi|x}^k\right) \le \Gamma_{\xi|x}^t \right\} \quad (12)$$

Let $L_{\xi|x}^k$ be a Cholesky factor of $\Sigma_{\xi|x}^k$, i.e., $\Sigma_{\xi|x}^k = L_{\xi|x}^k (L_{\xi|x}^k)^{\mathrm{T}}$. A finite set of unit direction vectors $\mathcal{V} := \{v_j\}_{j=1}^J \in \mathbb{R}^m$ with $\| v_j \|_2 = 1$ is chosen, where $J$ denotes the number of supporting halfspaces. In this letter, $J$ is set to 8. Subsequently, a polyhedral outer approximation of $\hat{\mathcal{W}}_{\text{CAUS}}^k(x)$ is constructed by the intersection of these supporting halfspaces, as follows:

$$\mathcal{W}_{\text{CAUS}}^{k,t}(x) = \left\{ \tilde{\xi}_t \mid v_j^{\mathrm{T}} (L_{\xi|x}^k)^{-1} \left(\tilde{\xi}_t - \boldsymbol{\mu}_{\xi|x}^k\right) \le \sqrt{\Gamma_{\xi|x}^t}, \forall j = 1, \ldots, J \right\} \quad (13)$$

Equivalently, $\mathcal{W}_{\text{CAUS}}^{k,t}(x)$ can be written as the standard nonempty and bounded polytope:

$$\mathcal{W}_{\text{CAUS}}^{k,t}(x) = \left\{ \tilde{\xi}_t \mid \boldsymbol{D}_k \tilde{\xi}_t \le \boldsymbol{d}_k \right\} \quad (14)$$

Finally, the CAUS is shown in Fig. 1(c), and can be formulated as:

$$\mathcal{W}_{\text{CAUS}} = \mathcal{W}_{\text{CAUS}}^1 \otimes \mathcal{W}_{\text{CAUS}}^2 \ldots \otimes \mathcal{W}_{\text{CAUS}}^T, \quad \mathcal{W}_{\text{CAUS}}^t = \bigcup_{k=1}^{K} \mathcal{W}_{\text{CAUS}}^{k,t} \quad (15)$$

**Proposition 1** *Assume that* $\{\tilde{\xi}_t \mid \tilde{x} = x\} \sim \mathbb{G}(\tilde{\xi}_t \mid \tilde{x} = x)$ *and* $\{\xi_t^{(i)}\}_{i=1}^{N_s}$ *utilized to compute* $\Gamma_{\xi|x}^t$ *are i.i.d. and drawn from the same conditional distribution. Then,* $\mathbb{P}(\tilde{\xi}_t \in \mathcal{W}_{\text{CAUS}}^t \mid \tilde{x} = x) \ge 1-\varepsilon$. *Moreover, for any decision s, if* $f_t(s,\xi) \le 0$ *for all* $\tilde{\xi}_t \in \mathcal{W}_{\text{CAUS}}^t$, *then* $\mathbb{P}\left(f_t(s,\tilde{\xi}_t) \le 0 \mid \tilde{x} = x\right) \ge 1-\varepsilon$.

**Proof.** Based on (12), by construction, $\xi_t \in \bigcup_{k=1}^K \hat{\mathcal{W}}_{\text{CAUS}}^{k,t}(x)$ if and only if $C_t(\xi;x) \le \Gamma_{\xi|x}^t$. Let $Z_i := C_t(\xi_t^{(i)};x)$ for $i$=1, 2, …, $N_s$ and $Z_{N_s+1} := C_t(\tilde{\xi};x)$. Then $Z_1, Z_2, \ldots, Z_{N_s+1}$ are iid conditional on $\tilde{x} = x$. With $\Gamma_{\xi|x}^t = Z_{(\kappa)}$, the $\kappa$-th order statistic of $Z_1, Z_2, \ldots, Z_{N_s}$, exchangeability implies $\mathbb{P}(Z_{N_s+1} \le Z_{(\kappa)} \mid \tilde{x} = x) \ge \kappa/(N_s+1) \ge 1-\varepsilon$. Hence $\mathbb{P}(\tilde{\xi} \in \bigcup_{k=1}^K \hat{\mathcal{W}}_{\text{CAUS}}^{k,t}(x) \mid \tilde{x} = x) \ge 1-\varepsilon$. Since $\hat{\mathcal{W}}_{\text{CAUS}}^{k,t}(x) \subseteq \mathcal{W}_{\text{CAUS}}^{k,t}(x)$, and therefore $\mathbb{P}(\tilde{\xi} \in \bigcup_{k=1}^K \mathcal{W}_{\text{CAUS}}^{k,t}(x) \mid \tilde{x} = x) \ge 1-\varepsilon$. Finally, if $f_t(s,\xi) \le 0$ for all $\tilde{\xi} \in \mathcal{W}_{\text{CAUS}}^t$, then $\{\tilde{\xi} \in \mathcal{W}_{\text{CAUS}}^t\} \subseteq \{f_t(s,\tilde{\xi}_t) \le 0\}$, implying $\mathbb{P}\left(f_t(s,\xi) \le 0 \mid \tilde{x} = x\right) \ge 1-\varepsilon$.■

As shown in (15), the handling of the uncertainty set represented as the UoS typically requires enumerating all subsets when searching for the worst-case realization. In a multi-period setting, this leads to $K^T$ candidate scenarios, which severely degrades computational efficiency. To address this issue, a mixed integer representation is proposed in (16):

$$\mathcal{W}_{\text{CAUS}} = \left\{ \tilde{\xi}_t \mid \alpha_{k,t} \in \{0,1\}, \sum_{k=1}^{K} \alpha_{k,t} = 1, \sum_{k=1}^{K} \alpha_{k,t} \left( D_k \tilde{\xi}_t - d_k \right) \le 0, \forall t \right\} \quad (16)$$

where $\alpha_{k,t} \in \{0,1\}$ are auxiliary variables. For any given subset, it can be uniquely identified by a binary indicator $\alpha_{k,t}$. Additionally, the bilinear terms $\alpha\tilde{\xi}$ can be reformulated as mixed linear constraints via Big-M formulation. The reformulated uncertainty set is then given in (17):

$$\begin{aligned} \mathcal{W}_{\text{CAUS}} = \Big\{ & \tilde{\xi}_t \mid \alpha_{k,t} \in \{0,1\}, \ \sum_{k=1}^{K} \alpha_{k,t} = 1, \\ & \sum_{k=1}^{K} \left( \boldsymbol{D}_k \boldsymbol{w}_{k,t} - \alpha_{k,t} \boldsymbol{d}_k \right) \le 0, \\ & -M\left(1-\alpha_{k,t}\right) \le \tilde{\xi}_t - \boldsymbol{w}_{k,t} \le M\left(1-\alpha_{k,t}\right), \\ & -M\alpha_{k,t} \le \boldsymbol{w}_{k,t} \le M\alpha_{k,t}, \forall t \Big\} \end{aligned} \quad (17)$$

where $\boldsymbol{w}_{k,t}$ represents auxiliary variables, and $M$ is a sufficiently large constant.

## III. Two-Stage Robust Unit Commitment

Building upon the proposed CAUS, its application to a two-stage robust unit commitment problem is demonstrated in this section, formulated as follows:

$$\begin{aligned} &\min_{X,Y,\tilde{\xi}} \left\{ \mathbf{c}^{\mathrm{T}} \mathbf{X} + \max_{\tilde{\xi}} \min_{Y} \mathbf{b}^{\mathrm{T}} \mathbf{Y} \right\} \\ &\text{s.t. } \mathbf{A}\mathbf{X} \le \mathbf{h} \\ &\quad \mathbf{E}\mathbf{X} + \mathbf{G}\mathbf{Y} \le \mathbf{l} - \mathbf{U}\tilde{\xi} \\ &\quad \tilde{\xi} \in \mathcal{W}_{\mathrm{CAUS}} \end{aligned} \tag{18}$$

where $\mathbf{X}$ denotes the first-stage decision variables, including the unit commitment variables, and $\mathbf{Y}$ denotes the second-stage recourse variables. $\mathbf{c}$ and $\mathbf{b}$ are the coefficient vectors within the objective, while $\mathbf{h}$ and $\mathbf{l}$ are coefficient vectors within the constraints. $\mathbf{A}$, $\mathbf{E}$, $\mathbf{G}$, and $\mathbf{U}$ are the constraint coefficient matrices. Additionally, $\tilde{\xi}$ is the renewable generation uncertainty vector, which belongs to the uncertainty set $\mathcal{W}_{\mathrm{CAS}}$.

The two-stage unit commitment problem can be decomposed into a master problem and a subproblem (SP). For an uncertainty set of the form presented in (15), each SP requires solving $K^T$ SPs to enumerate all subset combinations. In contrast, with the formulation given in (17), the number of decision variables is observed to grow only linearly with $K$ and $T$, which can improve computational efficiency.

## IV. Case Study

The effectiveness of the proposed CAUS is numerically validated utilizing a modified IEEE 118-bus system with three wind farms. Three methods, including RO, stochastic optimization (SO), and deterministic optimization (DO), are evaluated. For RO, the proposed CAUS is compared against the box set (BS), the UoS [7], and the partition-combine set (PCS) [2]. For SO, a sample average approximation is solved utilizing 1000 scenarios generated from a Gaussian distribution. All simulations are solved utilizing Gurobi 12.0. Additionally, 10000 out-of-sample realizations are generated by sampling forecast errors from the conditional distribution. The out-of-sample risk reliability is measured as the fraction of these realizations for which the recourse problem remains feasible under the fixed first-stage decisions. In this letter, $N_s$=10000, $\varepsilon$=0.05 are utilized.

The trade-off among cost, out-of-sample risk reliability, and computational time across DO, SO, and RO is summarized in Table I. The most optimistic strategy is provided by DO, yielding the lowest cost but only 53.14% out-of-sample risk reliability. Risk reliability is increased to 81.55% by SO, while its objective is about 1.69% higher than that of DO; however, a computational time of 9040.01s is required. For RO, 100% reliability is achieved by both BS and PCS, while 99.73% is achieved by UoS. Reductions of about 0.23% and 0.20% in the BS objective are realized by PCS and UoS, respectively. The BS objective is further reduced by 3.41% by CAUS, indicating that conservativeness is limited through the incorporation of side information. Additionally, CAUS is computationally competitive, requiring 301.5s, which represents a 15.50% reduction compared to PCS, a performance comparable to UoS, and nearly a 30-fold speedup relative to SO. As shown in Fig. 2, the results indicate that a controllable reliability guarantee is provided by the proposed method through the confidence level. Out-of-sample risk reliability is improved by increasing the confidence level. This shows that high out-of-sample reliability can be achieved by the proposed CAUS with limited additional conservativeness.

TABLE I

Comparative Results of Different Methods

| Method | | Cost ($10^6$\$) | Out-of-Sample Risk Reliability | Computational Time(s) |
|---|---|---|---|---|
| DO | | 1.7092 | 53.14% | 4.37 |
| SO | | 1.7381 | 81.55% | 9040.01 |
| RO | BS | 1.8729 | 100% | 21.17 |
| | PCS | 1.8686 | 100% | 356.8 |
| | UoS | 1.8692 | 99.73% | 286.4 |
| | **CAUS** | **1.8090** | **98.89%** | **301.5** |

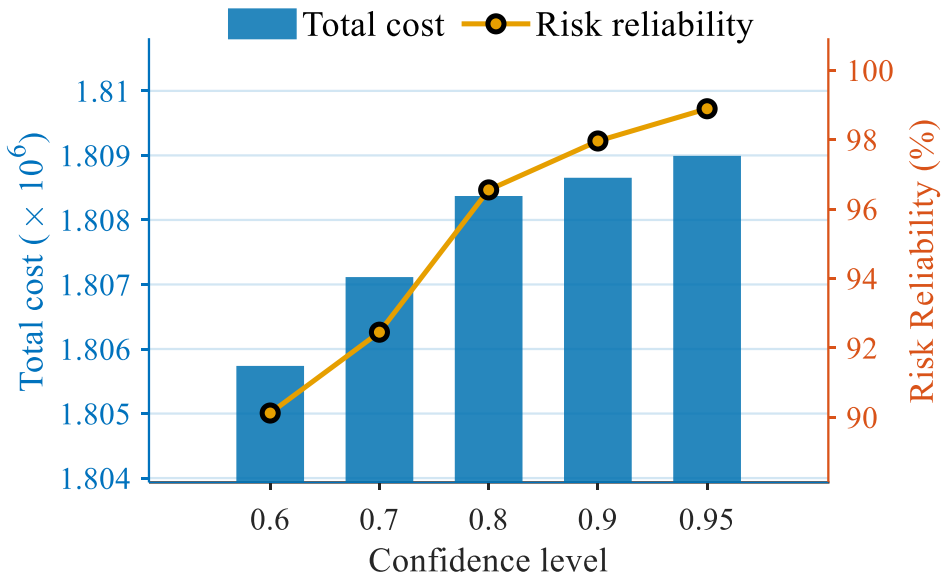

Fig. 2. Computational results of the proposed method.

## V. Conclusion

In this letter, a data-driven CAUS is proposed for application in robust dispatch of power systems. Side information is effectively incorporated by a conditional GMM to model irregular historical distributions, leveraging a UoS representation. To enhance scalability, a mixed integer linear reformulation is developed, enabling worst-case evaluation without explicit subset enumeration. Empirical evidence, obtained through case studies on robust unit commitment, demonstrates both reduced conservativeness and competitive computational performance for the proposed method.